\documentstyle[12pt]{article}
\input amssym.def
\input amssym
\begin{document}
\def \Z{\Bbb Z}
\def \C{\Bbb C}
\def \R{\Bbb R}
\def \Q{\Bbb Q}
\def \N{\Bbb N}
\def \wt{{\rm wt}}
\def \tr{{\rm tr}}
\def \span{{\rm span}}
\def \Res{{\rm Res}}
\def \Res{{\rm QRes}}
\def \End{{\rm End}}
\def \E{{\rm End}}
\def \Ind {{\rm Ind}}
\def \Irr {{\rm Irr}}
\def \Aut{{\rm Aut}}
\def \Hom{{\rm Hom}}
\def \mod{{\rm mod}}
\def \ann{{\rm Ann}}
\def \<{\langle} 
\def \>{\rangle} 
\def \t{\tau }
\def \a{\alpha }
\def \e{\epsilon }
\def \l{\lambda }
\def \L{\Lambda }
\def \g{\gamma}
\def \b{\beta }
\def \om{\omega }
\def \o{\omega }
\def \c{\chi}
\def \ch{\chi}
\def \cg{\chi_g}
\def \ag{\alpha_g}
\def \ah{\alpha_h}
\def \ph{\psi_h}
\def \be{\begin{equation}\label}
\def \ee{\end{equation}}
\def \bex{\begin{example}\label}
\def \eex{\end{example}}
\def \bl{\begin{lem}\label}
\def \el{\end{lem}}
\def \bt{\begin{thm}\label}
\def \et{\end{thm}}
\def \bp{\begin{prop}\label}
\def \ep{\end{prop}}
\def \br{\begin{rem}\label}
\def \er{\end{rem}}
\def \bc{\begin{coro}\label}
\def \ec{\end{coro}}
\def \bd{\begin{de}\label}
\def \ed{\end{de}}
\def \pf{{\bf Proof. }}
\def \voa{{vertex operator algebra}}

\newtheorem{thm}{Theorem}[section]
\newtheorem{prop}[thm]{Proposition}
\newtheorem{coro}[thm]{Corollary}
\newtheorem{conj}[thm]{Conjecture}
\newtheorem{exa}[thm]{Example}
\newtheorem{lem}[thm]{Lemma}
\newtheorem{rem}[thm]{Remark}
\newtheorem{de}[thm]{Definition}
\newtheorem{hy}[thm]{Hypothesis}
\makeatletter
\@addtoreset{equation}{section}
\def\theequation{\thesection.\arabic{equation}}
\makeatother
\makeatletter

\newcommand{\rw}{\rightarrow}
\newcommand{\n}{\:^{\circ}_{\circ}\:}

\begin{center}{\Large \bf Certain associative algebras similar to $U(sl_{2})$
and Zhu's algebra $A(V_{L})$}

\vspace{0.5cm}
Chongying Dong\footnote{Supported by NSF grant DMS-9303374 and a
research grant from the Committee on Research, UC Santa Cruz.},
Haisheng Li and Geoffrey Mason\footnote{Supported by NSF grant
DMS-9401272 and a research grant from the Committee on Research, UC
Santa Cruz.}\\ Department of Mathematics, University of California,
Santa Cruz, CA 95064
\end{center}

\begin{abstract}
It is proved that Zhu's algebra for vertex operator algebra
associated to a positive-definite even lattice of rank one is a 
finite-dimensional semiprimitive
quotient algebra of certain associative algebra introduced by Smith.
Zhu's algebra for vertex operator algebra associated to any
positive-definite even lattice is also calculated and is related to 
a generalization of Smith's algebra. 
\end{abstract}

\section{Introduction}
The recently developed vertex operator algebra theory encopes Lie
algebras (both finite-dimensional simple Lie algebras and
infinite-dimensional affine Lie algebras), groups (Lie groups and
finite groups), codes and lattices.  Many important algebras have
appeared to be closely related to certain vertex operator
algebras. For instance, the Griess algebra is a substructure of Frenkel,
Lepowsky and Muerman's Moonshine vertex operator algebra
$V^{\natural}$ whose full symmetry group is the Monster group.  Our
naive purpose in this paper is to relate a new interesting class of
associative algebras appeared in a different field to certain vertex
operator algebras.

In [S], Smith studied an interesting class of associative algebras $R(f)$ 
parameterized by a 
polynomial $f(x)$. Briefly, for any polynomial $f(x)$,
Smith defined an associative algebra $R(f)$ with three generators $A, B, H$ 
with defining relations $HA-AH=A, HF-FH=-F, AB-BA=f(H)$. As one of the main 
results in [S], it was shown that for some 
$f(x)$, $R(f)$ behaves much like $U(sl_{2})$ in terms of the complete 
reducibility of any finite-dimensional $R(f)$-modules. 
Following [BB], Hodges and Smith established an equivalence between the
category of $R(f)$-modules and the category of sheaves of left
$\cal{D}$-modules in [H] and [HS] where ${\cal D}$ is a sheaf of
rings on a siutable finite topological space.  

In this paper we find a relation between algebra $R(f)$ and vertex operator
algebras $V_L$ associated to positive definite even lattices $L$ of rank
1 by realizing Zhu algebras $A(V_L)$ as certain semisimple quotients
of $R(f)$ for certain $f.$ This relation suggests us to study
a generalization of $R(f)$ associated to any positive definite
even lattice which leads to the calculation of $A(V_L)$ in general.
In fact we characterize $A(V_L)$ by generators and realtions.
It is expected
that the realtion between generalization of $R(f)$ and $A(V_L)$
will help us to understand the vertex operator algebra structure 
of $V_L$ in terms of ${\cal D}$-modules on algebraic curves
(see [BD] and [HL]).

In [Z], an associative algebra
$A(V)$ was introduced for any
 vertex operator algebra $V$ so that there is a 1-1
correspondence between the set of equivalence classes of irreducible
$V$-modules and the set of equivalence classes of $A(V)$-modules. 
If $V$ is the irreducible
highest weight $\hat{\frak{g}}$-module of level $\ell$ with lowest weight 
$0$ where $\frak g$ is a finite-dimensional 
semisimple Lie algebra (cf. [FZ], [L1]), Zhu's algebra is isomorphic to $U({\frak{g}})$ if
 $\ell$ is generic. If $\ell$ is a positive integer, 
Zhu's algebra is isomorphic to a
semiprimitive quotient algebra of $U({\frak{g}})$ [FZ]. On the other hand, from FKS construction [FLM]
of basic modules for affine Lie algebras of types $A, D, E$, $L(1,0)=V_{L}$, where $L$ is
the root lattice of $\frak{g}$. So Zhu's
algebra $A(V_L)$ for an arbitrary  positive definite even lattice $L$
is very ``close'' to the universal
enveloping algebra of a semisimple Lie algebra in some sense. 
It is natural to expect certain ${\cal D}$-modules entering the picture
in the spirit of [BB].
 
Here is the precise structure of $A(V_L)$ if 
$L={\Z}\alpha$ be a rank-one lattice with $\<\alpha,\alpha\>=2k$
for some positive integer $k$. $A(V_{L})$ in this case is
a semiprimitive quotient ring of $R(f_{k})$ with
$$f_{k}(x)={2k\over (2k-1)!}x(4k^{2}x^{2}-1)(4k^{2}x^{2}-4)\cdots
(4k^{2}x^{2}-(k-1)^{2}).$$ 
In the case $\<\a,\a\>=2$ $A(V_L)$ has been computed previously in [Lu].

The paper is organized as follows. In Section 2 we study the algebra
$R(f_k)$ and prove that a certain quotient $\bar R_k$ of $R(f_k)$
is a semi-simple algebra whose all irreducible representations
are given explicitly. A generalization $\bar A(L)$ of $R(f)$ 
associated to an positive definite lattice $L$ is also investigated.
We establish the semisimplicity of $\bar A(L)$ and construct all
its irrducible modules. In Section 3 we prove that
$\bar R_k$ and $A(V_{\Z\a})$ are isomorphic if $\<\a,\a\>=2k$
and that $\bar A(L)$ and $A(V_L)$ are isomorphic for arbitrary 
positive definite even lattice $L.$    

\section{Certain associative  algebras similar to $U(sl_{2})$}

Let $g(x)\in {\C}[x]$ be any polynomial in $x$. The associative algebra
$R(g)$ [S] is generated by 
$\{A, B, H\}$ subject to relations:
$$HA-AH=A,\;\; HB-BH=-B,\;\; AB-BA=g(H).$$
In fact, $R(g)$ is a ${\Z}$-graded algebra with 
$\deg A=1=-\deg B$ and $\deg H=0,$ and 
$A(g)$ has a basis $\{ B^{m}H^{n}A^{k}|  m,n,k\in {\Z}_{+}\}$. 

Let $P$ be the subalgebra of $R(g)$ generated by $A$ and $H$. Then 
$AP=\{Aa|a\in P\}$ is a two-sided ideal of $P$ and $P=AP\oplus {\C}[H]$. 
For any complex
number $\lambda$, let ${\C}v_{\lambda}$ be the $1$-dimensional
$P$-module such that $APv_{\lambda}=0, Hv_{\lambda}=\lambda
v_{\lambda}$.  Define the Verma
$R(g)$-module as follows:
$$V(\lambda)=R(g)\otimes_P{\C}v_{\lambda}.$$ Then $V(\lambda)$ has a
unique maximal proper submodule.  Denote by $L(\lambda)$ the
irreducible quotient module of $V(\lambda)$.

Let $u(x)$ be a polynomial of degree $n+1$ and take  
$g={1\over 2}(u(x+1)-u(x))$. Similar to $U(sl_{2})$, $R(g)$ has a central element
$$\Omega=AB+BA+{1\over 2}\left(u(H+1)+u(H)\right)$$
such that the center of $R(g)$ is isomorphic to ${\C}[\Omega]$.
Furthermore, $\Omega$ acts on $V(\lambda)$ as a scalar $u(\lambda+1)$.

The following two propositions can be found in [S].
\bp{ns}
(a) The finite-dimensional simple $R(g)$-modules are precisely the
modules $L(\lambda)=V(\lambda)/B^{j}v_{\lambda}$ where $j\in {\N}$ is
minimal such that $u(\lambda+1)=u(\lambda+1-j)$.

(b)The number of simple modules of dimension $j$ equals
$$|\{\lambda\in {\C}|u(\lambda+1)=u(\lambda+1-j), \mbox{ and j is the
least such element of }{\N}\}|$$ which is less than or equal to 
$\deg u-1=\deg g.$
\ep

\bp{psemi}
Suppose that for each $j\in {\N}$, there are precisely $\deg g$
simple modules of dimension $j$. Then every finite-dimensional
$R(g)$-module is semisimple.
\ep

Set $h_j(x)=g(x)+g(x-1)+\cdots +g(x-j)$ for $j\in \Z_+.$ Observe that
$h_j(x)=\frac{1}{2}(u(x)-u(x-j)).$  The next corollary
is implicit in [S]. 
\bc{csemi1} If each root of $h_j(x)=0$ for $j\in\Z_+$ has
multiplicity $1$ and any two equations $h_i(x)=0, h_j(x)=0$
do not have common solutions, then every
finite-dimensional $R(g)$-module is semisimple.
\ec

{\bf Proof.} Since $\deg h_{j-1}(x)=\deg (u(x+1)-u(x+1-j))=\deg g(x),$
$u(x+1)=u(x+1-j)$ has exactly $\deg g$ solutions for fixed $j\in \N.$
By the assumption, $u(x+1)=u(x+1-j)$ and $u(x+1)=u(x+1-i)$ have no
common solutions if $i\ne j.$ Thus
for each solution $\lambda$ of $u(x+1)=u(x+1-j),$ $j$ is the minumum
such that $u(\lambda+1)=u(\lambda+1-j).$ By Proposition \ref{ns}
$R(g)$ has exactly $\deg g$ simple modules of dimension $j$. 
Use Proposition \ref{psemi} to finish the proof.$\;\;\;\;\Box$

Similar to the complete reducibility of an integrable module for a 
semisimple Lie algebra (cf. [K]) we have

\bc{csemi2} If each root of $h_j(x)=0$ for $j\in\Z_+$ has
multiplicity $1$ and any two equations $h_i(x)=0, h_j(x)=0$
do not have common solutions, then every
$R(g)$-module on which $A$ and $B$ are locally nilpotent is
semisimple.
\ec

{\bf Proof.} We first prove that any nonzero $R(g)$-module $M$ on which 
$A$ and $B$ are locally nilpotent contains a 
finite-dimensional simple $R(g)$-module $V$.
 
Set $M^{0}=\{u\in M|Au=0\}$. It is clear that $HM^{0}\subseteq
M^{0}$.  Since $A$ is locally nilpotent on $M$, $M^{0}\ne 0$. Let
$0\ne v\in M^{0}$.  Then there exists a nonnegative integer $r$ such
that $B^{r}u\ne 0$ and $B^{r+1}v=0$.  Since
\begin{eqnarray*}
& &2^{r+1}A^{r+1}B^{r+1}\nonumber\\
&=&(\Omega-u(H))(\Omega-u(H-1))\cdots (\Omega-
u(H-r))\nonumber\\
&=&(2BA+u(H+1)-u(H))(2BA+u(H+1)-u(H-1))\cdots \nonumber\\
& &\cdots (2BA+u(H+1)-u(H-r))
\end{eqnarray*}
(cf. [DLM3] or [S]), we obtain $$(u(H+1)-u(H))(u(H+1)-u(H-1))\cdots (u(H+1)-u(H-r))v=0.$$
Consequently, there is a $\lambda\in {\C}$ and a $0\ne u\in M^{0}$
such that $Hu=\lambda u$.  Then $u$ generates a highest weight
$R(g)$-module $V$. By Corollary \ref{csemi1}, $V$ is simple. 

Now let $W$ be the sum of all
finite-dimensional simple $R(g)$-submodules of $M$. Then it suffices to prove
$M=W$.  Suppose $W\ne M$. Then $\bar{M}=M/W$ is a nonzero
$R(g)$-module which has a finite-dimensional simple
submodule $M^{1}/W$ where $M^{1}$ is a submodule of
$M$. Let $u\in M^{1}$ such that $u+W$ is a nonzero highest weight
vector. Then $R(g)u$ is contained in a finite-dimensional
$R(g)$-submodule $L$ of $M^{1}$. It follows from Corollary
\ref{csemi1} that $L$ is a direct sum of finite-dimensional simple
$R(g)$-modules so that $u$ is contained in $W$. It is a
contradiction. The proof is complete.$\;\;\;\;\Box$

Smith [S] gave an instructive example for $g(x)=(x+1)^{n+1}-x^{n+1}$.
Motivated by vertex operator algebras associated to positive definite even
lattices of rank one we consider the associative
algebras $R_k=R(g_k)$ for any positive integer $k$ where
$$g_{k}(x)=\frac{1}{(2k-1)!}2kx(4k^{2}x^{2}-1)(4k^{2}x^{2}-4)\cdots
(4k^{2}x^{2}-(k-1)^{2}).$$ Then $R_k$
is generated by $\{A, B, H\}$ subject to relations:
\begin{eqnarray*}
& &HA-AH=A,\;\; HB-BH=-B,\\
& &AB-BA=
\frac{1}{(2k-1)!}2kH(4k^{2}H^{2}-1)(4k^{2}H^{2}-4)\cdots (4k^{2}H^{2}-(k-1)^{2}).
\end{eqnarray*}

Notice that $R_{1}=U(sl_{2})$. If $k=2$, the corresponding polynomial
$u_{2}(x)={16\over 3}\left(x-{1\over 2}\right)^{4} -{10\over
3}\left(x-{1\over 2}\right)^{2}.$ For any positive integer $k$ it may be 
possible to prove that any
finite-dimensional $R_{k}$-module is completely reducible. 
But in this paper we only prove the result for $k=2.$

\bl{k2}
Every $R_{2}$-module $M$ on which $A$ and $B$ are locally nilpotent is
semisimple.
\el

{\bf Proof.} Since
 $g_{2}(x)={2\over 3}x(16x^{2}-1)={32\over 3}x^{3}-{2\over 3}x$,
for any $r\in {\Z}_{+}$, we have
\begin{eqnarray*}
& &g_{2}(x)+g_{2}(x-1)+\cdots +g_{2}(x-r)\\
&=&{32\over 3}\left(x^{3}+(x-1)^{3}+\cdots +(x-r)^{3}\right)
-{2\over 3}\left(x+(x-1)+\cdots +(x-r)\right)\\
&=&{(r+1)\over 3}(2x-r)\left(16x^{2}-16rx+8r^{2}+8r-1\right).
\end{eqnarray*}
Notice that $16x^{2}-16rx+8r^{2}+8r-1=0$ has two distinct
noreal solutions.  If $\lambda$ satisfies two equations
$16x^{2}-16rx+8r^{2}+8r-1=0$ and $16x^{2}-16sx+8s^{2}+8s-1=0$ for
$r>s>0$, then $\lambda$ satisfies the difference equation
$16(r-s)x-8(r^{2}-s^{2})-8(r-s)=0$.  Consequently, $\lambda$ is
rational. It is a contradiction. Thus any two of the equations
$g_{2}(x)+g_{2}(x-1)+\cdots +g_{2}(x-j)=0$ for $j=0,1,\cdots$ do not
have common solutions. The result  follows from Corollary \ref{csemi2}
immediately.$\;\;\;\;\Box$

For a fixed positive integer $k$, we define $\bar{R}_{k}$ to be the
quotient algebra of $R_{k}$ modulo the two-sided ideal generated by
$(1-2H)A$.  Next, we shall prove that $\bar{R}_k$ is a semisimple
algebra. We need the following lemma.

\bl{lone}
Let $M$ be any $\bar{R}_{k}$-module. Then $H$ is locally finite on $M$.
\el

{\bf Proof.} We first show that $M$ contains a
nonzero $\bar{R}_{k}$-submodule on which $H$ is semisimple.
If $AM=0$, then $g_{k}(H)M=0$ because $[A,B]=g_{k}(H)$. Since the
multiplicity of any root of $g_{k}(x)=0$ is one, $H$ is semisimple on
$M$. If $AM\ne 0$, let $0\ne u=Av$ for some $v\in M.$
 Then $(1-2H)u=(1-2H)Av=0$. Thus,
$Hu={1\over 2}u$. Then $\bar{R}_{k}u$ is a nonzero
$\bar{R}_{k}$-submodule on which $H$ semisimple.

Let $F(M)$ be the maximal $\bar{R}_{k}$-submodule (of $M$) on which
$H$ is locally finite.  If $F(M)=M$, we are done. Otherwise, consider
the quotient module $W=M/F(M)$. Let $W_1$ be a nonzero 
submodule of $W$ on which $H$ is locally finite. Write
$W_{1}=M_{1}/F(M)$ where $M_{1}\supset F(M)$ is a submodule of $M.$ 
Then $H$ is locally finite on $M_{1}$. This is a
contradiction. $\;\;\;\;\Box$

\bt{tone}
$\bar{R}_{k}$ is  semisimple (finite-dimensional) and all inequivalent
irreducible $\bar{R}_{k}$-modules are $L({n\over 2k})$ for 
$n\in {\Z}, -(k-1)\le n\le k$.
\et

{\bf Proof.} To prove that $\bar{R}_{k}$ is semisimple is 
equivalent to prove that any $\bar{R}_{k}$-module is semisimple.  Let
$M$ be any $\bar{R}_{k}$-module. Define $$M_{\lambda}=\{u\in M|
(H-\lambda)^{n}u=0\;\;\mbox{for some positive integer }n.\}$$ Then by
Lemma \ref{lone} we have $M=\oplus_{\lambda\in {\C}}M_{\lambda}$. It
is easy to see that $AM_{\lambda}\subseteq M_{\lambda+1}$ and
$BM_{\lambda}\subseteq M_{\lambda-1}.$ Let $u\in M_{\lambda}$ for some
$\lambda\in {\C}$. Use relation $(1-2H)Au=0$ to get $HAu={1\over 2}Au.$ 
Thus $Au\in M_{{1\over 2}}$ and $Au=0$ if $\lambda\ne -{1\over 2}$.  Set
$M^{1}=\oplus_{\lambda\in {1\over 2}+{\Z}}M_{\lambda}$ and
$M^{2}=\oplus_{\lambda\notin {1\over 2}+{\Z}}M_{\lambda}$. Then
$M^{1}$ and $M^{2}$ are two submodules of $M$ and $M=M^{1}\oplus
M^{2}$. Clearly $AM^{2}=0$. From the proof of
Lemma \ref{lone} we know that all the eigenvalues of $H$ in 
$M^{2}$ are roots of $g_k(x).$ That is,
$$M^2=\oplus_{n=-k+1}^kM_{\frac{n}{2k}}.$$
As $1+\frac{n}{2k}$ is not in $\{\frac{m}{2k}|-k+1\leq m\leq k\}$
for any $-k+1\leq n\leq k$ we see that $BM^2=0.$ So each $M_{\frac{n}{2k}}$ 
is a direct sum of 1-dimensional
irreducible $\bar{R}_{k}$-modules isomorphic to $L(\frac{n}{2k}).$
In fact $L(\frac{n}{2k})$ for $-k+1\leq n\leq k$ give all the inequivalent 
1-dimensional irreducible $\bar R_k$-modules. 

It remains  to show that $M^{1}$ is semisimple.
Take $\lambda\in {1\over 2}+{\Z}$ with $\lambda\ne \pm {1\over 2}$. Then
$\lambda-1\ne -{1\over 2}$ and $AM_{\lambda}=0=ABM_{\lambda}$. This yields
$g_{k}(H)M_{\lambda}=0$. If $M_{\lambda}\ne 0$, we must have
$\lambda={n\over 2k}$ for some $-(k-1)\le n\le k-1$. This is 
a contradiction.  Thus
$M_{\lambda}=0$ if $\lambda\in {1\over 2}+{\Z}, \lambda\ne \pm {1\over
2}$ and $M^{1}=M_{-{1\over 2}}\oplus M_{{1\over
2}}$. Consequently, $AM_{{1\over 2}}=0=BM_{-{1\over 2}}$. If $u\in
M_{-{1\over 2}}$, we get $g_{k}(-{1\over
2})u=g_{k}(H)u=[A,B]u=BAu$. Since $g_{k}(-{1\over 2})\ne 0$, $A$ is an
injective map from $M_{-{1\over 2}}$ to $M_{{1\over 2}}$ and $B$ is a
surjective map from $M_{{1\over 2}}$ to $M_{-{1\over 2}}$. Similarly,
$A$ is a surjective map from $M_{-{1\over 2}}$ to $M_{{1\over 2}}$ and
$B$ is an injective map. Since $H$ acts on $AM$ as a scalar ${1\over
2}$, $H$ acts on $M_{{1\over 2}}\;(=AM_{-{1\over 2}})$ as a scalar
${1\over 2}$ and $H$ acts on $M_{-{1\over 2}}\;(=BM_{{1\over 2}})$ as
a scalar $-{1\over 2}$. Let $u_{i}$ for $i\in I$ be a basis of
$M_{{1\over 2}}$.  Then $M^{1}=\oplus _{i\in
I}({\C}u_{i}+{\C}Bu_{i})$. It is easy to see that for any $i\in I$,
$({\C}u_{i}+{\C}Bu_{i})$ is a submodule which is isomorphic to
$L({1\over 2})$.  Thus $M^{1}$ is a direct sum of several copies of
2-dimensional irreducible $\bar{R}_{k}$-module $L({1\over
2})$. Clearly, $L(\frac{1}{2})$ is the only irreducible $\bar R_k$-module
whose dimension is greater than 1. $\;\;\;\;\Box$

\br{rweight}
{}From the proof of Theorem \ref{tone}, we see that $H$ is semisimple
 on any $\bar{R}_{k}$-module $M$ such that any $H$-weight
$\lambda$ satisfies $2k\lambda\in {\Z}, -k\le 2k\lambda \le k$.
\er

Let $x$ be an indeterminant. Then we define ${x\choose 0}=1$ and
${x\choose r}={1\over r!}x(x-1)\cdots (x-r+1)$ for any positive
integer $r$.  Then for any $m\ge n\in {\Z}_{+}$, we have:
\begin{eqnarray}\label{2.6}
\sum_{i=0}^{n}{n\choose i}{x\choose m-i}={x+n\choose m}.
\end{eqnarray}

Let $L$ be any positive-definite even lattice and 
Let $\hat{L}$ be the canonical central extension of $L$ by the cyclic
group $\< \pm 1\>$:
\begin{eqnarray}\label{2.7}
1\;\rightarrow \< \pm 1\>\;\rightarrow \hat{L}\;\bar{\rightarrow}
L\;\rightarrow 1 
\end{eqnarray}
with the commutator map $c(\alpha,\beta)=(-1)^{\< \alpha,\beta\>}$ for
$\alpha,\beta \in L$. Let $e: \hat L\to L$ be a section such that
$e_0=1$ and $\epsilon: L\times L\to \<\pm 1\>$ be
the corresponding 2-cocycle. Then $\epsilon(\a,\b)\epsilon(\b,\a)=(-1)^{\<\a,\b\>},$
\begin{equation}\label{2c}
\e(\a,\b)\e(\a+\b,\gamma)=\e(\b,\gamma)\e(\a,\b+\gamma)
\end{equation} 
and $e_{\a}e_{\b}=
\e(\a,\b)e_{\a+\b}$ for $\a,\b,\gamma\in L.$ 

Set ${\bf h}={\C}\otimes_{{\Z}}L$.  For any $\alpha,\beta\in L$, we define
$g_{\alpha,\beta}(x)=0$ if $\<\alpha,\beta\>\ge 0$, and define
\begin{eqnarray}\label{gdef}
g_{\alpha,\beta}(x)=\sum_{r=0}^{{\<\alpha,\alpha\>\over 2}-1}
{{\<\alpha,\alpha\>\over 2}-1\choose r}{x\choose -\<\alpha,\beta\>-1-r}
={x+{\<\alpha,\alpha\>\over 2}-1\choose -\<\alpha,\beta\>-1}
\end{eqnarray}
if $\<\alpha,\beta\>< 0$.  

We define an associative algebra
$A(L)$ generated by $E_{\alpha}$ $(\alpha\in L)$ and ${\bf h}$ subject to
relations:
\begin{eqnarray}
& &E_{0}=1 \mbox{  (the identity)};\label{edef1}\\
& &hh'-h'h=0\;\;\;\mbox{ for any }h,h'\in {\bf h};\label{edef2}\\
& &hE_{\alpha}-E_{\alpha}h=\<h,\alpha\>E_{\alpha};\label{edef3}
\end{eqnarray}
for any $\alpha,\beta\in L$. 

\br{rlie}
Let $L$ be the root lattice of a semisimple Lie algebra $\frak{g}$
with the set $\Delta$ of roots.  Then the subalgebra of $A(L)$
generated by ${\bf h}$ and $E_{\alpha}$ for $\alpha\in \Delta$ modulo
$$E_{\alpha}E_{\beta}-E_{\beta}E_{\alpha}=
E_{\alpha+\beta}g_{\alpha,\beta}(\alpha)\e(\a,\b)$$
is isomorphic to the universal enveloping algebra $U(\frak{g})$ of
$\frak{g}$.
\er

Next we define a quotient algebra $\bar{A}(L)$ of $A(L)$ modulo the
following relations:
\begin{eqnarray}
& &\left(\alpha-{\<\alpha,\alpha\>\over 2}\right)E_{\alpha}=0;\label{e1}\\
& &E_{\alpha}E_{\beta}=0\;\;\;\mbox{if }\<\alpha,\beta\> > 0;\label{e2}\\
& &E_{\alpha}E_{\beta}
=E_{\alpha+\beta}\e(\a,\b){\alpha+{\<\alpha,\alpha\>\over 2}\choose -\<\alpha,\beta\>}
\;\;\;\mbox{if }\<\alpha,\beta\>\le 0\label{e3}
\end{eqnarray}
for $\alpha,\beta\in L$.

\bp{pal}
Any $\bar{A}(L)$-module is completely reducible. That is,
$\bar{A}(L)$ is a (finite-dimensional) semisimple algebra.
\ep

{\bf Proof.} For any $\alpha\in L$, define $\bar{A}_{\alpha}(L)$ to be
the subalgebra of $\bar{A}(L)$ generated by $\alpha, E_{\alpha},
E_{-\alpha}$. Note that $\e(\a,-\a)\e(-\a,\a)=1$ and 
$${\a+\frac{\<\a,\a\>}{2}\choose\<\a,\a\>}-{-\a+\frac{\<\a,\a\>}{2}\choose\<\a,\a\>}={\a+\frac{\<\a,\a\>}{2}-1\choose\<\a,\a\>-1}.$$
Set $2k=\<\alpha,\alpha\>.$
Then $\bar{A}_{\alpha}(L)$ is isomorphic to a
quotient algebra of $\bar{R}_{k}$ by sending $A,B, H$ to 
$\e(\a,-\a)E_{\alpha},E_{-\alpha}, {1\over 2k}\alpha$,
respectively.  Since any $\bar{A}(L)$-module $M$ is a
direct sum of irreducible $\bar{A}_{\alpha}(L)$-modules, $\alpha$ is semisimple
on $M.$ Note that  $L$ spans ${\bf h}$. So ${\bf h}$ is semisimple on $M.$ 
Denote by $M_{\lambda}$ the ${\bf h}$-weight space of weight $\lambda$.
Then $M=\oplus_{\lambda\in {\bf h}}M_{\lambda}.$
By Remark \ref{rweight} if $M_{\lambda}\ne 0$ then 
$|\<\lambda,\alpha\>|\le {1\over
2}\<\alpha,\alpha\>$ and $\<\lambda,\alpha\>\in\Q$ 
for any $\alpha\in L$.  

For any $u\in M_{\lambda}$, we set $M(u)=\oplus_{\alpha\in
L}{\C}E_{\alpha}u$. Note that  $E_{0}=1$. It follows from the relations 
(\ref{edef3}),
(\ref{e2}) and (\ref{e3}) that $M(u)$ is a submodule containing $u.$
Relation
(\ref{e1}) gives $\alpha E_{\alpha}u={1\over
2}\<\alpha,\alpha\>E_{\alpha}u$ for $\alpha\in L.$ Use (\ref{edef3})
to obtain 
$$\alpha
E_{\alpha}u=(\<\alpha,\lambda\>+\<\alpha,\alpha\>)E_{\alpha}u.$$ 
This gives $(\<\alpha,\lambda\>+{1\over 2}\<\alpha,\alpha\>)E_{\alpha}u=0$.
Then either $\<\alpha,\lambda\>+{1\over 2}\<\alpha,\alpha\>=0$ or
$E_{\alpha}u=0$.  If $\<\alpha,\lambda\>+{1\over 2}
\<\alpha,\alpha\>\ne 0,$ $E_{\alpha}u$ must be 0. Since $L$ is
positive-definite, there are only finitely many $\alpha\in L$
satisfying the relation $\<\alpha,\lambda\>+{1\over 2}
\<\alpha,\alpha\>=0$. Thus $M(u)$ is finite-dimensional. Suppose
$E_{\alpha}u\ne 0$. Then $\<\alpha,\lambda\>+{1\over 2}
\<\alpha,\alpha\>=0$ and
$E_{-\alpha}E_{\alpha}u=\e(-\a,\a)u$ (see
definition (\ref{e3})).  This proves that $M(u)$ is a
finite-dimensional irreducible $\bar{A}(L)$-module.  As a result,
$M$ is a direct sum of finite-dimensional irreducible
$\bar{A}(L)$-modules. This completes the proof.$\;\;\;\;\Box$

\br{rcon}
Let $M$ be an $\bar{A}(L)$-module. Then it follows from the proof
of Proposition \ref{pal} and Remark \ref{rlie} 
 that any ${\bf h}$-weight $\lambda$ of $M$ is in 
the dual lattice $L^{\circ}=\{\lambda\in{\bf h}|\<L,\lambda\>\subset\Z\}$
 of $L$ and satisfies the relation:
$|\<\lambda,\alpha\>|\le {1\over 2}\<\alpha,\alpha\>$ for any
$\alpha\in L$, which is equivalent to the relation:
$\<\lambda+\alpha,\lambda+\alpha\>\ge \<\lambda,\lambda\>$ for any
$\alpha\in L$.
\er

Next, we construct all irreducible $\bar{A}(L)$-modules. Set
$$S=\{\lambda\in L^{\circ}|\<\lambda+\alpha,\lambda+\alpha\>\ge
\<\lambda,\lambda\>\;\;
\mbox{ for any }\alpha\in L\}.$$
For any $\lambda\in S$, we define $$\Delta(\lambda)=\{\alpha\in
L|\<\lambda+\alpha,\lambda+\alpha\> = \<\lambda,\lambda\>\}.$$ 
Then
for any $\lambda\in S, \alpha\in \Delta(\lambda)$, we have
$\lambda+\alpha\in S$ because
$$\<\lambda+\alpha+\beta,\lambda+\alpha+\beta\>\geq \<\lambda,\lambda\>
=\<\lambda+\alpha,\lambda+\alpha\>$$ for any $\beta\in L$. 

For any $\lambda\in S$, let $M^{\lambda}$ be a vector space with a
basis $\{u_{\alpha}^{\lambda}|\alpha\in
\Delta(\lambda)\}$.
Define
$$hu_{\alpha}^{\lambda}=\<\lambda+\alpha,h\>u_{\alpha}^{\lambda}\;\;\;\mbox{ for }
h\in {\bf h}, 
\alpha\in \Delta(\lambda).$$ 
For $\beta\in L, \alpha\in \Delta(\lambda)$,
we define $E_{\beta}u_{\alpha}^{\lambda}=\e(\a,\b)
u_{\alpha+\beta}^{\lambda}$ if 
 $\beta \in \Delta(\alpha+\lambda)$ and $E_{\beta}u_{\alpha}^{\lambda}=0$ 
otherwise.

\bp{pm} The following hold:

(1) The vector space $M^{\lambda}$ together with the defined action is
an irreducible $\bar{A}(L)$-module for any $\lambda\in S.$

(2) $M^{\lambda_{1}}\simeq
M^{\lambda_{2}}$ if and only if $\lambda_{2}=\lambda_{1}+\alpha$ for
some $\alpha\in \Delta(\lambda_{1})$. 

(3) Any irreducible $\bar{A}(L)$-module is isomorphic to $M^{\lambda}$ for some
$\lambda\in S$.
\ep

{\bf Proof.} For (1) we need to establish the defining relations
(\ref{edef1})-(\ref{e3}) of $\bar A(L).$  Relations (\ref{edef1})-(\ref{edef3})
are clear. Let $\alpha\in\Delta(\l)$ and $\b\in L.$ If $\b\in \Delta(\l+\a)$
then $2\<\l+\a,\b\>+\<\b,\b\>=0.$ Thus
$$(\b-\frac{\<\b,\b\>}{2})E_{\b}u_{\a}^{\l}
=(\b-\frac{\<\b,\b\>}{2})\e(\b,\a)u_{\a+\b}^{\l}=(\<\l+\a+\b,\b\>-\frac{\<\b,\b\>}{2})\e(\b,\a)u_{\a+\b}^{\l}=0.$$
If  $\b\not\in \Delta(\l+\a),$ $E_{\b}u_{\a}^{\l}=0.$ This gives 
relation (\ref{e1}). 

We now show realtion (\ref{e2}).
Let $\beta_{1},\beta_{2}\in L$ such that 
$\<\b_1,\b_2\>>0.$ We have to show that $E_{b_1}E_{b_2}u_{\a}^{\l}=0$ for any
$\alpha\in \Delta(\lambda),$ or equivalently, either  $\beta_2\not
\in \Delta(\a+\l)$ or $\beta_1\not\in \Delta(\a+\b_1+\l).$ 
If $\beta_2\in \Delta(\a+\l)$ and  $\beta_1\in \Delta(\a+\b_2+\l)$ then 
$$\<\lambda+\alpha+\b_1+\beta_2,\lambda+\alpha+\b_1+\beta_2\>=\<\lambda+\alpha+\beta_2,\lambda+\alpha+\beta_2\>
=\<\lambda+\alpha,\lambda+\alpha\>=\<\lambda,\lambda\>$$ 
and $\b_1\in\Delta(\a+\b_2+\l).$ 
This implies that 
$$\<\lambda+\alpha,\beta_{1}\>+{1\over 2}\<\beta_{1},\beta_{1}\>+\<\beta_{1},\beta_{2}\>
=0.$$
Since $\<\lambda+\alpha,\beta_{1}\>+{1\over 2}\<\beta_{1},\beta_{1}\>\ge 0$
we see that  $\<\beta_{1},\beta_{2}\> \leq 0$. This is a contradiction.

It remains to show (\ref{e3}). Let $\beta_{1},\beta_{2}\in L$ such that 
$\<\b_1,\b_2\>\leq 0.$ Take $\a\in\Delta(\l).$ There are 3 cases:

Case 1: $\beta_2,\beta_1+\beta_2\in \Delta(\a+\l).$ From the previous 
paragraph, $\b_1\in \Delta(\a+\b_2+\l).$ Use (\ref{2c}) to obtain
\begin{eqnarray*}
& &E_{\beta_{1}}E_{\beta_{2}}u_{\alpha}^{\lambda}=\e(\b_1,\b_2+\a)\e(\b_2,\a)u_{\beta_{1}+\beta_{2}+\alpha}^{\lambda}\\
& &\ \ \ =\e(\b_1,\b_2)\e(\b_1+\b_2,\a)u_{\beta_{1}+\beta_{2}+\alpha}^{\lambda}\\
& &\ \ \ =\e(\b_1,\b_2)E_{\beta_1+\beta_2}u_{\alpha}^{\lambda}.
\end{eqnarray*}
Also note that 
$\<\beta_{1},\beta_{2}\> \le 0$. So ${{\<\b_1,\a+\l\>+\frac{\<\b_1,\b_1\>}{2}}
\choose -\<\b_1,\b_2\>}=1.$ (\ref{e3}) is true in this case.

Case 2: $\beta_1+\b_2\not\in \Delta(\a+\l).$
Then $E_{\beta_{1}+\beta_{2}}u_{\alpha}^{\lambda}=0.$ From the discussion 
before we see that either  $\beta_2\not
\in \Delta(\a+\l)$ or $\beta_1\not\in \Delta(\a+\b_1+\l).$ So 
$E_{\beta_{1}}E_{\beta_{2}}u_{\alpha}^{\lambda}=0.$

Case 3: $\beta_2\not\in \Delta(\a+\l)$ and $\beta_1+\b_2\in \Delta(\a+\l).$
In this case $E_{\beta_{1}}E_{\beta_{2}}u_{\alpha}^{\lambda}=0,$
$\<\lambda+\alpha,\beta_{2}\>+{1\over 2}\<\beta_{2},\beta_{2}\> >0$
and 
$$\<\lambda+\alpha,\beta_{1}+\beta_{2}\>
+{1\over 2}\<\beta_{1}+\beta_{2},\beta_{1}+\beta_{2}\> =0.$$
we have 
$$0\le \<\lambda+\alpha,\beta_{1}\>+{1\over 2}\<\beta_{1},\beta_{1}\>
<-\<\beta_{1},\beta_{2}\>.$$
Then 
${\<\lambda+\alpha,\beta_{1}\>+{1\over 2}\<\beta_{1},\beta_{1}\>\choose 
-\<\beta_{1},\beta_{2}\>}=0,$ as desired.

It is clear that $M^{\lambda}$ is irreducible.

(2) follows from the definition of $M^{\l}$ and 
(3) follows from the proof of Proposition \ref{pal} immediately.
$\;\;\;\;\Box$

\bc{cnumber}
There is a one-to-one correspondence between set of equivalence
classes of irreducible $\bar{A}(L)$-modules and the set of cosets of
$L$ in $L^{\circ}$.
\ec

{\bf Proof.} For any $\lambda_{1},\lambda_{2}\in S$, we define
$\lambda_{1}\equiv \lambda_{2}$ if $\lambda_{2}=\lambda_{1}+\alpha$
for some $\alpha \in \Delta(\lambda_{1})$. If
$\lambda_{1},\lambda_{2}\in S$ satisfies
$\lambda_{2}-\lambda_{1}=\alpha\in L$, then it follows from the
definitions of $S$ and $\Delta(\lambda_{1})$ that $\alpha\in
\Delta(\lambda_{1})$.  Then we have an equivalent relation on $S$,
which is the restriction of the congruence relation of $L^{\circ}$ modulo
$L$. On the other hand, for any coset $\lambda+L\in L^{\circ}/L$, since $L$
is positive-definite, there is an element $\beta\in \lambda+L$ (as a
set) with a minimal norm. Clearly $\beta\in S$.  The corollary
follows from Proposition \ref{pm}.  $\;\;\;\;\Box$

\newpage

\section{Zhu's algebra $A(V_{L})$}

First we recall from [FLM] the explicit construction of vertex operator algebra
$V_{L}.$ Let $L$ be a
positive-definite even lattice. Set ${\bf h}={\C}\otimes_{{\Z}}L$ and
extend the ${\Z}$-form on $L$ to ${\bf h}$.  Let $\hat{{\bf h}}={\bf
C}[t,t^{-1}]\otimes {\bf h} \oplus {\bf C}c$ be the affinization of
${\bf h}$, i.e., $\hat{{\bf h}}$ is a Lie algebra with commutator
relations:
\begin{eqnarray*}
[t^{m}\otimes h,t^{n}\otimes h']=m\delta _{m+n,0}\< h,h'\> c\;\;\mbox{for
}h,h'\in {\bf h};m,n\in {\Z};
\end{eqnarray*}
\begin{eqnarray*}
[\hat{{\bf h}},c]=0.
\end{eqnarray*}
We also use the notation $h(n)=t^{n}\otimes h$ for $h\in {\bf h}, n\in
{\Z}$.

Set 
\begin{eqnarray*}
\hat{{\bf h}}^{+}=t{\bf C}[t]\otimes {\bf h};\;\;\hat{{\bf
h}}^{-}=t^{-1}{\bf C}[t^{-1}]\otimes {\bf h}.
\end{eqnarray*} 
Then
$\hat{{\bf h}}^{+}$ and $\hat{{\bf h}}^{-}$ are abelian subalgebras of
$\hat{{\bf h}}$.  
Let $U(\hat{{\bf h}}^{-})=S(\hat{{\bf h}}^{-})$ be
the universal enveloping algebra of $\hat{{\bf h}}^{-}$. Consider the
induced $\hat{{\bf h}}$-module
\begin{eqnarray*}
M(1)=U(\hat{{\bf h}})\otimes _{U({\bf C}[t]\otimes {\bf h}\oplus {\bf C}c)}{\bf
C}\simeq S(\hat{{\bf h}}^{-})\;\;\mbox{(linearly)},\end{eqnarray*}
where ${\bf C}[t]\otimes {\bf h}$ acts trivially on ${\bf C}$ and $c$ acts
on ${\bf C}$ as multiplication by 1.

Recall from (\ref{2.7}) that
$\hat{L}$ be the canonical central extension of $L$ by the cyclic
group $\< \pm 1\>.$ Form the induced $\hat{L}$-module
\begin{eqnarray*}
{\C}\{L\}={\bf C}[\hat{L}]\otimes _{\< \pm 1\>}{\bf C}\simeq
{\bf C}[L]\;\;\mbox{(linearly)},\end{eqnarray*}
where ${\bf C}[\cdot]$ denotes the group algebra and $-1$ acts on ${\bf
C}$ as multiplication by $-1$. For $a\in \hat{L}$, write $\iota (a)$ for
$a\otimes 1$ in ${\bf C}\{L\}$. Then the action of $\hat{L}$ on ${\C}
\{L\}$ is given by: $a\cdot \iota (b)=\iota (ab)$ and $(-1)\cdot \iota
(b)=-\iota (b)$ for $a,b\in \hat{L}$.

Furthermore we define an action of ${\bf h}$ on ${\bf C}\{L\}$ by:
$h\cdot \iota (a)=\< h,\bar{a}\> \iota (a)$ for $h\in {\bf h},a\in
\hat{L}$. Define $z^{h}\cdot \iota (a)=z^{\< h,\bar{a}\> }\iota (a)$.

The untwisted space associated with $L$ is defined to be
\begin{eqnarray*}
V_{L}={\bf C}\{L\}\otimes _{{\bf C}}M(1)\simeq {\bf C}[L]\otimes
S(\hat{{\bf h}}^{-})\;\;\mbox{(linearly)}.\end{eqnarray*}
Then $\hat{L},\hat{{\bf h}},z^{h}\;(h\in {\bf h})$ act naturally on
$V_{L}$ by
acting on either ${\bf C}\{L\}$ or $M(1)$ as indicated above.

For $h\in {\bf h}$ set
\begin{eqnarray*}
h(z)=\sum _{n\in {\Z}}h(n)z^{-n-1}.
\end{eqnarray*}
We use a normal ordering procedure, indicated by open colons, which
signify that in the enclosed expression, all creation operators $h(n)$
$(n<0)$,$a\in \hat{L}$ are to be placed to the left of all
annihilation operators $h(n),z^{h}\;(h\in {\bf h},n\ge 0)$. For $a \in
\hat{L}$, set
\begin{eqnarray*}
Y(\iota (a),z)=\n e^{\int (\bar{a}(z)-\bar{a}(0)z^{-1})}az^{\bar{a}}\n.
\end{eqnarray*}
Let $a\in \hat{L};\;h_{1},\cdots,h_{k}\in {\bf h};n_{1},\cdots,n_{k}\in {\Z}\;(n_{i}> 0)$. Set
\begin{eqnarray*}
v= \iota (a)\otimes h_{1}(-n_{1})\cdots h_{k}(-n_{k})\in V_{L}.\end{eqnarray*}
Define vertex operator 
\begin{eqnarray*}
Y(v,z)=\n \left({1\over (n_{1}-1)!}({d\over
dz})^{n_{1}-1}h_{1}(z)\right)\cdots \left({1\over (n_{k}-1)!}({d\over
dz})^{n_{k}-1}h_{k}(z)\right)Y(\iota (a),z)\n 
\end{eqnarray*}
This gives us a well-defined linear map
\begin{eqnarray*}
Y(\cdot,z):& &V_{L}\rightarrow (\mbox{End}V_{L})[[z,z^{-1}]]\nonumber\\ 
& &v\mapsto Y(v,z)=\sum _{n\in {\Z}}v_{n}z^{-n-1},\;(v_{n}\in {\rm
End}V_{L}).
\end{eqnarray*}

Let $\{\;\a_{i}\;|\;i=1,\cdots,d\}$ be an orthonormal basis of ${\bf
h}$ and set
\begin{eqnarray*}
\omega ={1\over 2}\sum _{i=1}^{d}\a_{i}(-1)\a_{i}(-1)\in V_{L}.
\end{eqnarray*}
Then $Y(\omega,z)=\sum_{n\in {\bf Z}}L(n)z^{-n-2}$ gives rise to a
representation of the Virasoro algebra on $V_{L}$ and 
\begin{eqnarray}
& &L(0)\left(\iota(a)\otimes h_{1}(-n_{1})\cdots
h_{n}(-n_{k})\right)\nonumber \\
&=&\left({1\over 2}\< \bar{a},\bar{a}\>+n_{1}+\cdots+n_{k}\right)
\left(\iota(a)\otimes h_{1}(-n_{1})\cdots
h_{k}(-n_{k})\right).
\end{eqnarray}

The following theorem was due to Borcherds [B] and Frenkel, Lepowsky
and Meurman [FLM].

\bt{tBFLM} 
$(V_{L},Y,{\bf 1},\omega)$ is a vertex operator algebra. 
\et

Recall the dual lattice $L^{\circ}$ of $L$ from Section 2.  Let
$\{\beta_{i}+L|i=1,2,\cdots, n\}$ be a full set of representives of
cosets of $L$ in $L^{\circ}$.  Then $V_{\beta_{i}+L}$ is an
irreducible $V_{L}$-module for each $i$ [FLM] and it was also proved in [D]
that these are all irreducible $V_{L}$-modules up to equivalence;
see also [DLM1].

Define the Schur polynomials $p_{r}(x_{1},x_{2},\cdots)$ $(r\in
{\Z}_{+})$ in variables $x_{1},x_{2},\cdots$ by the following
equation:
\begin{eqnarray}\label{eschurd}
\exp \left(\sum_{n= 1}^{\infty}\frac{x_{n}}{n}y^{n}\right)
=\sum_{r=0}^{\infty}p_{r}(x_1,x_2,\cdots)y^{r}.
\end{eqnarray}
For any monomial $x_{1}^{n_{1}}x_{2}^{n_{2}}\cdots x_{r}^{n_{r}}$ we
have an element $h(-1)^{n_{1}}h(-2)^{n_{2}}\cdots
h(-r)^{n_{r}}{\bf 1}$ in $V_{L}$ for $h\in{\bf h}.$ Then for any polynomial
$f(x_{1},x_{2}, \cdots)$, $f(h(-1), h(-2),\cdots){\bf 1}$ is a
well-defined element in $V_{L}$.  In particular,
$p_{r}(h(-1),h(-2),\cdots){\bf 1}$ for $r\in {\Z}_{+}$ are elements
of $V_{L}$.

Suppose $a,b\in \hat{L}$ such that $\bar{a}=\alpha,\bar{b}=\beta$. 
Then
\begin{eqnarray}\label{erelation}
Y(\iota(a),z)\iota(b)&=&z^{\<\alpha,\beta\>}\exp\left(\sum_{n=1}^{\infty}
\frac{\alpha(-n)}{n}z^{n}\right)\iota(ab)\nonumber\\
&=&\sum_{r=0}^{\infty}p_{r}(\alpha(-1),\alpha(-2),\cdots)\iota(ab)
z^{r+\<\alpha,\beta\>}.
\end{eqnarray}
Thus 
\begin{eqnarray}\label{eab1}
\iota(a)_{i}\iota(b)=0\;\;\;\mbox{ for }i\ge -\<\alpha,\beta\>.
\end{eqnarray}
Especially, if $\<\alpha,\beta\>\ge 0$, we have
$\iota(a)_{i}\iota(b)=0$ for all $i\in {\Z}_{+}$, and if
$\<\alpha,\beta\>=-n<0$, we get
\begin{eqnarray}\label{eab}
\iota(a)_{i-1}\iota(b)=p_{n-i}(\alpha(-1),\alpha(-2),\cdots)\iota(ab)
\;\;\;\mbox{ for }i\in {\Z}_{+}.
\end{eqnarray}
It is well-known that $V_{L}$ is generated by $u_{n}$ for 
$u\in \{\iota(a), \alpha(-1)|a\in \hat{L}\}, n\in {\Z}$.

In [Z], an associative algebra $A(V)$ was introduced for any vertex
operator algebra $V=\oplus_{n\in\Z}V_n$ such that there is a 1-1
correspondence between the set of equivalence classes of irreducible
$V$-modules (without assuming that homogeneous subspaces are
finite-dimensional) and the set of all equivalence classes of
irreducible $A(V)$-modules on which the central element (obtained
{}from the Virasoro element $\omega$) as a scalar.  

Here is the
definition of $A(V).$ Define two bilinear products $*$ and $\circ$ on
$V$ as follows:
\begin{eqnarray*}
& &u*v={\rm Res}_{z}\frac{(1+z)^{n}}{z}Y(u,z)v=\sum_{i=0}^{\infty}{n\choose i}
u_{i-1}v;\\
& &u\circ v={\rm Res}_{z}\frac{(1+z)^{n}}{z^{2}}Y(u,z)v
=\sum_{i=0}^{\infty}{n\choose i}u_{i-2}v
\end{eqnarray*}
for any $u\in V_{n}, v\in V$.  Denote by $O(V)$ the linear span of all
$u\circ v$ for $u,v\in V$ and set $A(V)=V/O(V)$. Then $A(V)$ is an 
associative algebra under $*$ with identity ${\bf 1}+O(V).$
For any (weak) $V$-module $M$, we define
[DLM2] $$\Omega(M)=\{u\in M|a_{m}u=0\;\;\mbox{ for any }a\in V_{n},
m>n-1\}.$$ Then $\Omega(M)$ is a natural $A(V)$-module with $a+O(V)$
acting as $o(a)=a_{n-1}$ for $a\in V_{n}$. This gives is a 1-1
correspondence between the set of equivalence classes of irreducible
$V$-modules (without assuming that homogeneous subspaces are
finite-dimensional) and the set of all equivalence classes of
irreducible $A(V)$-modules [Z].

The following relations from [Z] 
will be useful later:
\begin{eqnarray}
& &u*v\equiv {\rm Res}_{z}\frac{(1+z)^{{\wt v}-1}}{z}Y(v,z)u
\equiv\sum_{i=0}^{\infty}{{\wt v}-1\choose i}v_{i-1}v\;\;\mbox{mod }O(V);\label{tra}\\
& &u*v-v*u\equiv {\rm Res}_{z}(1+z)^{{\wt u}-1}Y(u,z)v
\equiv \sum_{i=0}^{\infty}{{\wt u}-1\choose i}u_{i}v\;\;\mbox{mod }O(V).\label{ez2}
\end{eqnarray}

\bt{tlattice}
Let $L={\Z}\alpha$ be a one-dimensional lattice with $|\alpha|^{2}=2k$
for $k\in {\N}$ and let $V_{L}$ be the vertex operator algebra
associated with $L$. Then Zhu's algebra $A(V_{L})$ is isomorphic to
$\bar{R}_{k}=R_{k}/\< (1-2H)A\>$.
\et

{\bf Proof.} In this case $\hat L$ is a direct product of $L$ and
$\<\pm 1\>$ and we regard $L$ as a subgroup of $\hat L.$ The section
$e$ is the identity map and the cocycle $\e(\cdot,\cdot)$ is trivial. 
Let $a\in \hat L$ such that $\bar a=\alpha.$ 
Set $e=\iota(a), f=\iota(a^{-1}),
h=\alpha(-1){\bf 1}$.  {}From the construction of $V_{L}$ and
(\ref{eab}) we have
\begin{eqnarray*}
& &h_{i}e=2k\delta_{i,0}e,\;\; h_{i}f=-2k\delta_{i,0}f,\; h_{i}h=2k \delta_{i,1}{\bf 1};\\ 
& &e_{i}f=p_{2k-1-i}(h_{-1},h_{-2},\cdots){\bf 1}
\end{eqnarray*}
for any $i\in {\Z}_{+}$. Since $\wt h=1, \wt e=\wt f=k$, by (\ref{ez2}) we have
\begin{eqnarray*}
& &h*e-e*h\equiv {\rm Res}_{z}{1\over z}Y(h,z)e\equiv h_{0}e\equiv 2ke
\;\;\mbox{mod }O(V_{L});\\
& &h*f-f*h\equiv {\rm Res}_{z}{1\over z}Y(h,z)f\equiv h_{0}f\equiv -2kf
\;\;\mbox{mod }O(V_{L});\\
& &e*f-f*e\equiv \sum_{i=0}^{k-1}{k-1\choose i}e_{i}f\nonumber\\
& &\hspace{2cm}\equiv \sum_{i=0}^{k-1}{k-1\choose i}p_{2k-1-i}(h(-1),h(-2),\cdots)
{\bf 1}
\;\;\;\mbox{mod }O(V_{L}).
\end{eqnarray*}

For any  $r\ge 1$ and $u\in V_{L}$ we have
$$(h(-r-1)+h(-r))u={\rm Res}_{z}\frac{(1+z)^{\wt h}}{z^{r+1}}u\in O(V)$$
(cf. [FZ]).  
Thus $h(-r-1)\equiv -h(-r)u$  mod $O(V)$ for any $u\in V_{L}$. Then
for any $r\in {\Z}_{+}, n_{1},\cdots, n_{r}\in {\N}$, we have:
\begin{eqnarray}\label{extra1}
h(-n_{1})\cdots h(-n_{r}){\bf 1}+O(V_{L})=(-1)^{n_{1}+\cdots +n_{r}+r}h^{r}+O(V_{L}).
\end{eqnarray}
Let $\bar{p}_{r}(x)=p_{r}(x,-x,x,\cdots)$, {\it i.e.,} substitute
$x_{n}$ by $(-1)^{n-1}x$ for any $n\ge 1$. Since 
$$\exp \left(\sum_{n= 1}^{\infty}\frac{(-1)^{n-1}x}{n}y^{n}\right)
=(1+y)^x=\sum_{r\geq 0}{x\choose r}y^r$$
we see that $\bar{p}_{r}(x)={x\choose r}={1\over r!}x(x-1)\cdots
(x+1-r)$  for $r\in {\Z}_{+}.$
Hence
\begin{eqnarray*}
& &e*f-f*e\equiv \sum_{i=0}^{k-1}{k-1\choose i}\bar{p}_{2k-1-i}(h)\equiv 
\sum_{i=0}^{k-1}{k-1\choose i}{h\choose 2k-1-i}\nonumber\\
& &\equiv {h+k-1\choose 2k-1}\;\;\mbox{mod }O(V_{L}).
\end{eqnarray*}
Here ${h+k-1\choose 2k-1}$ is understood to be $\frac{1}{(2k-1)!}(h+k-1)*(h+k-2)*\cdots *(h-k+1)$  and this consideration also applies in the next theorem. 
This gives an algebra homomorphism $\psi$ from $R_{k}$ to
$A(V_{L})$ such that $\psi(A)=e+O(V_{L}),\; \psi(B)=f+O(V_{L}),\;
\psi(H)={1\over 2k}h+O(V_{L})$.

Set $W_{\pm}=\sum_{n=0}^{\infty}e^{\pm n\alpha}\otimes M(1)$. Then
$V_{L}=W_{-}+W_{+}$.  Consider $G_{+}=\sum_{m\in {\Z}}{\C}e_{m}$ and
$G_{-}=\sum_{m\in {\Z}}{\C}f_{m}$. Then $G_{\pm}$ are abelian Lie
subalgberas of $G_{\pm}+\hat{{\bf h}}$ which is a Lie subalgebra
of $\End V_L.$  It is clear
that $W_{\pm}$ is generated by $G_{\pm}+\hat{{\bf h}}$ from the vacuum
vector, {i.e.,} $W_{\pm}=U(G_{\pm}+\hat{{\bf h}}){\bf 1}$. By PBW
theorem we have $W_{\pm}=U(G_{\pm}+\hat{{\bf h}}){\bf 1}=U(\hat{{\bf
h}}^{-})U(G_{\pm}){\bf 1}$.  Since $e_{m}{\bf 1}=0$ for $m\in
{\Z}_{+}$, $W_{+}$ can be generated from the vacuum ${\bf 1}$ by the
following operators:
\begin{eqnarray*}
& &{\rm Res}_{z}\frac{(1+z)^{k}}{z^{n+1}}Y(e,z)=\sum_{i=0}^{k}{k\choose i}e_{-n-1+i};\\
& &{\rm Res}_{z}\frac{(1+z)^{1}}{z^{n+1}}Y(h,z)=h(-n-1)+h(-n)
\end{eqnarray*}
for $n\in {\Z}_{+}$. Similarly for $W_-.$ This implies
that the map $\psi$ is onto.

A straightforward calculation gives $L(-1)\iota
(a)=\a(-1)\iota (a)$. Then 
\begin{equation}\label{extra}
L(-1)\iota
(a)=(\alpha(-1)+\alpha(0))\iota (a)-\alpha(0)\iota (a)
=(\alpha(-1)+\alpha(0))\iota (\alpha)-2k\iota (a). 
\end{equation}
Since $L(0)+L(-1)$ maps $V$ to $O(V)$ (see [Z]) and $h*\iota(a)=
(\a(-1)+\a(0))\iota(a)$ 
we have $-k\iota(a)=-L(0)\iota(a)\equiv L(-1)v=h*\iota(a)-2k\iota(a).$
This gives rise to the relation $-k\psi(A)=2k\psi(H)*\psi(A)-2k\psi(A)$ in $A(V_{L})$. That is,
$(1-2\psi(H))*\psi(A)=0$. This shows that $A(V_{L})$ is a quotient algebra of
$\bar{R}_{k}$.  Since $V_{L}$ already has irreducible modules
$V_{L+{n\over 2k}\alpha}$ for $-(k-1)\le n\le k$ and $\bar{R}_{k}$ has
exactly $2k$ irreducible modules, $A(V_{L})$ must be isomorphic to
$\bar{R}_{k}$.$\;\;\;\;\Box$

\br{rverma} It is interesting to know if one can construct a vertex operator 
algebra whose Zhu's algebra is exactly $R_{k}$. Notice that $\{e_{n},
f_{n}, h_{n}| n\in {\Z}\}$ generates a topological associative algebra
$A$ because of the infinite sum relation (\ref{erelation}), which is
not a linear Lie algebra. In order to construct such a vertex operator
algebra one may need to develop the notion of
generalized Verma modules for $A$ and establish some results on
generalized Verma modules.  Of course, if we took these for granted,
then we could have a vertex operator algebra with $R_{k}$ as its Zhu's
algebra.
\er

\bt{tvl}
Let $L$ be any positive-definite even lattice. Then
Zhu's algebra $A(V_{L})$ is isomorphic to $\bar{A}(L)$.
\et

{\bf Proof.} The proof is similar to that of Theorem \ref{tlattice}.
First, we establish an algebra homomorphism $\psi$ from $\bar{A}(L)$
onto $A(V_{L})$. Recall from Section 2 the section $e: L\to \hat L.$  
Define a linear map $\psi$ from $\bar{A}(L)$ to
$V_{L}$ as follows:
\begin{eqnarray*}
\psi(E_{\alpha})=\iota(e_{\alpha}), \psi(h)=h_{-1}{\bf 1}=h(-1)
\end{eqnarray*}
for $\alpha\in L, h\in {\bf h}$. From (\ref{ez2})  we have
$$h(-1)*h'(-1)-h'(-1)*h(-1)\equiv {\rm Res}_{z}\frac{1}{z}Y(h(-1),z)h'(-1)+O(V)=0$$
for $h,h'\in{\bf h}.$ 
Similarly, for $a\in \hat L, h\in {\bf h}$ with $\bar a=\alpha$
$$h(-1)*\iota(a)-\iota(a)*h(-1)\equiv h_{0}\iota(a)
\equiv \<h,\alpha\>\iota(a)\;\;\mbox{ mod }O(V_{L}).$$
The same calculation in (\ref{extra}) gives
$$L(-1)\iota(a)=\alpha(-1)*\iota(a)-\<\a,\a\>\iota(a).$$
Thus $L(-1)\iota(a)+L(0)\iota(a)=\alpha(-1)*\iota(a)-\frac{\<\a,\a\>}{2}\iota(a)\in O(V_L).$ This shows that relations (\ref{edef1})-(\ref{e1}) hold.

Let $a,b\in\hat L$ such that  $\bar a=\alpha$ and $\bar b=\beta.$ Then 
$$\iota(a)*\iota(b)=\sum_{i=0}^{m}{m-1\choose i}\iota(a)_{i-1}\iota(b)$$
where $m=\frac{\<\alpha,\alpha\>}{2}$ is the weight of $\iota(a)$.
If $\<\alpha,\beta\>> 0$, by (\ref{eab1}) we have
 $\iota(\alpha)_{i}\iota(\beta)=0$ for all $i\geq -1$, so that
$\iota(\alpha)*\iota(\beta)=0$ and (\ref{e2}) holds.
 
If $\<\alpha,\beta\>=-n\leq 0$, by (\ref{eab}) we have
$$\iota(a)*\iota(b)=\sum_{i=0}^{m}{m\choose i}p_{n-i}(\alpha(-1),\alpha(-2),\cdots)\iota(ab).$$ 
Now take $a=e_{\a}$ and $b=e_{\b}.$ 
As in (\ref{extra1}) we obtain
\begin{eqnarray*}
& &\iota(a)*\iota(b)+O(V_{L})\nonumber\\
&=&\sum_{i=0}^{m}{m\choose i}p_{n-i}(\alpha(-1),-\alpha(-1),\cdots)\iota(ab)+O(V_{L}).
\end{eqnarray*}
Since $\wt h(-1)=1$, by (\ref{tra})  we get 
$u*h(-1)+O(V_{L})=h(-1)u+O(V_{L})$ for any $u\in V_L.$ 
Then $h(-1)^k\iota(a)+O(V_L)=u*h(-1)^k$ for any nonnegative integer $k.$
As in the previous theorem 
 $h(-1)^k$ is understood to be the multiplication $*$ $k$ times.
Thus
\begin{eqnarray*}
& &\iota(a)*\iota(b)+O(V_{L})\nonumber\\
\nonumber\\
&=&\sum_{i=0}^{m}{m\choose i}p_{n-i}(\alpha(-1),-\alpha(-1),\cdots)\iota(ab)+O(V_{L})
\nonumber\\
&=&\iota(ab)*\left(\sum_{i=0}^{m}{m\choose i}p_{n-i}(\alpha(-1),-\alpha(-1),\cdots)\right)
+O(V_{L})\nonumber\\
&=&\iota(ab)*\left(\sum_{i=0}^{m}{m\choose i}{\alpha(-1)\choose n-i}\right)
+O(V_{L})\nonumber\\
&=&\iota(ab)*{\alpha(-1)+m\choose n}+O(V_{L}).
\end{eqnarray*}
So the relation (\ref{e3}) is true
and $\psi$ is an algebra homomorphism from $\bar{A}(L)$ into $A(V_{L})$.

Recall that $\{\a_1,...,\a_d\}$ is an orthonomal basis for ${\bf h}.$
Let $u=p(\a_{i}(-j))\iota(a)$ for 
$a\in \hat{L}$ and $p(x_{i,j})\in {\C}[x_{i,j}|1\le i\le d,j=1,2,\cdots]$. 
Then from the previous paragraph we see that 
\begin{eqnarray}
u\equiv p((-1)^{j-1}\a_{i}(-1))\iota(a)\equiv \iota(a)* (p(\a_{i}(-1)){\bf 1})
\;\;\mbox{mod }O(V).
\end{eqnarray}
Since such $u$ span $V_L$ (by the construction of $V_{L}$), 
$\psi$ is onto. Since $V_L$ already has $|L^{\circ}/L|$ modules (see
[FLM] and [D]), it follows from Proposition \ref{pal} that $\psi$ 
is an isomorphism.
$\;\;\;\;\Box$

\br{rdong}
Since there is a 1-1 correspondence between the set of equivalence classes of
irreducible $V_{L}$-modules and the set of equivalence classes of irreducible
$A(V_{L})$-modules, our result on $A(V_{L})$ gives an alternative approach to 
the classification of irreducible $V_{L}$-modules obtained in [D].
\er

\end{document}